\documentclass[12pt]{amsart}
\usepackage{somedefs,tracefnt,latexsym,,rawfonts,epsfig,amsfonts,amssymb,latexsy
m,enumerate,amscd}
\def\cal{\mathcal}
\def\Bbb{\mathbb}

\def\r{\rangle}
\def\l{\langle}
\def\t{\times}

\newtheorem{prop}{Proposition}[section]
\newtheorem{thm}{Theorem}[section]
\newtheorem{exm}{Example}[section]

\newtheorem{lemma}{Lemma}[section]
\newtheorem{cor}{Corollary}[section]

\newtheorem*{fjic}{(Fibered) Isomorphism conjecture}
\newtheorem{defn}{Definition}[section]

\newtheorem{rem}{Remark}[section]
\begin{document}
\title[Isomorphism conjecture] {On the isomorphism conjecture for groups 
acting on trees}
\author[S.K. Roushon]{S.K. Roushon}
\address{School of Mathematics\\
Tata Institute\\
Homi Bhabha Road\\
Mumbai 400005, India}
\email{roushon@math.tifr.res.in} 
\urladdr{http://www.math.tifr.res.in/\~\ roushon/}
\begin{abstract} We study the Fibered Isomorphism 
conjecture of Farrell and Jones for groups acting on trees. We show that 
under certain conditions the conjecture is true for groups acting on  
trees when the stabilizers satisfy the conjecture. These 
conditions are satisfied in several cases of the  
conjecture. We prove some general results on the conjecture for
the pseudoisotopy theory for groups acting on 
trees with residually
finite vertex stabilizers. In particular, we study situations when the 
stabilizers belong to the following classes of
groups: polycyclic groups, finitely generated nilpotent groups, closed
surface groups, finitely
generated abelian groups and virtually cyclic groups. 

Finally, we provide 
explicit examples of groups on which the results of this article 
can be applied and 
show that these groups were not considered before. Furthermore,  
we deduce that these groups are neither hyperbolic nor $CAT(0)$.\end{abstract}

\keywords{pseudoisotopy theory, $K$-theory, fibered isomorphism conjecture,
equivariant homology 
theory, group action on trees, $CAT(0)$-groups.}

\subjclass[2000]{Primary: 19D35, 57N37. Secondary: 55N91.}
\maketitle
\tableofcontents
\section{Introduction and statements of Theorems}
In this article we study the Fibered Isomorphism conjecture 
of Farrell and Jones for groups acting on trees. Originally the
conjecture was stated for the pseudoisotopy theory, algebraic $K$-theory
and for $L$-theory. 
Here we prove some results for the pseudoisotopy theory case
of the conjecture. The methods of proofs of our results hold for
the other theories under certain conditions. We make these
 conditions explicit here. It is well known that the pseudoisotopy 
version of the 
 conjecture yields computations of the Whitehead groups and lower 
$K$-groups of 
the associated groups and also implies the conjecture in 
lower $K$-theory (that is in dimension $\leq 1$). 
The main problem we are concerned with is the following.

\noindent
{\bf Problem.} {\it Assume the Fibered Isomorphism conjecture for the
 stabilizers of the action of a group on a tree. Show that the group also
satisfies the conjecture.}

A particular case of the problem is to show that the conjecture is true 
for a generalized free product $G_1*_HG_2$ ($HNN$-extension $G*_H$) of groups 
assuming   
it is true for the groups $G_1$, $G_2$ and $H$ ($G$ and $H$).

In [\cite{R1}, Reduction Theorem] we solved the problem when the 
edge stabilizers are trivial.
In this article we show it when the edge stabilizers are finite and in
addition either the vertex stabilizers are residually finite or 
the group surjects onto another group with certain properties. 
We also solve the
problem, under a certain condition, assuming that 
the vertex stabilizers are polycyclic.
This condition is both necessary and sufficient 
for the fundamental 
group to be subgroup separable when the vertex stabilizers are finitely
generated 
nilpotent. We further consider cases when the stabilizers are virtually cyclic or 
surface groups. Next, we prove a result when the stabilizers are finitely 
generated abelian. Graph of groups of the last type were studied before 
in \cite{S} where the 
$L$-theory was computed when the stabilizers are finitely generated free abelian
groups. We note that in this case the result in \cite{S} 
implies the $L$-theory version of the Isomorphism conjecture. 

Furthermore, a positive answer to the problem will imply that the Fibered 
Isomorphism conjecture is true for one-relator groups for all the three 
theories.

Before we come to the statements of the main results let us 
recall that a {\it graph of groups} consists of a graph 
$\cal G$ (that is an one-dimensional $CW$-complex) and 
to each vertex $v$ or edge $e$ of $\cal G$ there is associated a group 
${\cal G}_v$ (called the {\it vertex group} of the vertex $v$) 
or ${\cal G}_e$ (called 
the {\it edge group} of the edge $e$) respectively with the assumption that 
for each edge $e$ and for its two end vertices $v$ and $w$ (say)
there are injective group homomorphisms ${\cal G}_e\to {\cal G}_v$ and 
${\cal G}_e\to {\cal G}_w$. (If $v=w$ then also we demand two 
injective homomorphisms as above). The {\it fundamental group} 
$\pi_1({\cal G})$ of the graph $\cal G$ can be defined so that 
in the simple cases of graphs of groups where the graph has 
two vertices and one edge or one vertex and one edge the 
fundamental group is the amalgamated free product or the 
$HNN$-extension respectively. 
See \cite{DD} for some more on this subject.

Throughout the article we make the following conventions: 

1. A graph is assumed to be locally finite and connected. 

2. We use the same notation for a graph of groups and
its underlying graph. $V_{\cal G}$ and $E_{\cal G}$ denotes 
respectively the set of all vertices and edges of $\cal G$.   

3. And a group is assumed to be countable.

Recall that given two groups $G$ and $H$ the {\it wreath product} 
$G\wr H$ is by definition the semidirect product $G^H\rtimes H$ where the 
action of $H$ on $G^H$ is the regular action and $G^H$ denotes the 
direct sum of copies of $G$ indexed by the elements of $H$.

Now we state our main results.

In the statements of the results we say that the FIC$^P$ 
(FICwF$^P$) is true for
a group 
$G$ if the Farrell-Jones Fibered Isomorphism conjecture for the pseudoisotopy 
theory is true for the group $G$ (for $G\wr H$ for any finite group $H$). We 
mention here that the advantage of taking the FICwF$^P$ is that it passes to 
finite index over-groups (see Proposition \ref{product}).

\begin{thm} \label{residually} Let $\cal G$ be a graph of groups so that
the edge groups are finite. Then $\pi_1
({\cal G})$ satisfies the FICwF$^P$  
if one of the following conditions 
is satisfied.

(1). The vertex groups are residually finite and satisfy the FICwF$^P$. 

(2). There is a homomorphism                  
$f:\pi_1({\cal G})\to Q$ onto another group $Q$ so that the restriction
of $f$ to any vertex group has finite kernel and $Q$ satisfies the FICwF$^P$. 
\end{thm}

This Theorem is a special case of Proposition \ref{residually-prop} $(1)$ 
respectively Proposition \ref{theorem} $(b)$. See Section 7.

\begin{rem} {\rm There is a large class of residually
finite groups for which the FICwF$^P$  
is true. For example, any group which
contains a
group from the following examples as a subgroup of finite index is a
residually finite group satisfying the FICwF$^P$.

1. Polycyclic groups (\cite{FJ}).

2. Artin full braid groups (\cite{FR}, \cite{R3}).

3. Compact $3$-manifold groups (\cite{R1}, \cite{R2}).

4. Compact surface groups (\cite{FJ}).

5. Fundamental groups of hyperbolic Riemannian manifolds (\cite{FJ}).

6. Crystallographic groups (\cite{FJ}).}\end{rem}

Here we remark that by an email, dated July 11, 2008,
Wolfgang L\"{u}ck informed the author that the proof of the FIC$^P$ for a certain
class of polycyclic
groups as stated in [\cite{FJ}, Theorem 4.8] is not yet complete. However,
this result and the corresponding result in the $L$-theory case of the
conjecture were already used by several authors in published papers. The 
proof in the $L$-theory case is completed 
recently in \cite{BFL} and the pseudoisotopy case is in progress. 
We note here that
if a polycyclic group $G$ is also virtually nilpotent then the FIC$^P$ for $G$
is already proved in [\cite{FJ}, Lemma 4.1]

\begin{defn}\label{def1.1}
{\rm A graph of groups $\cal G$ is 
said
to satisfy the {\it intersection property} if for each 
connected subgraph of groups
${\cal G}'$ of $\cal G$, $\cap_{e\in E_{{\cal G}'}}{\cal G}'_e$
contains a subgroup which is normal in
 $\pi_1({\cal G}')$ and is of finite index in some edge group. We say  
$\cal G$ is of {\it finite-type} if the graph is finite and 
all the vertex groups are finite.}\end{defn}

\begin{thm} \label{ip} Let $\cal G$ be a graph of groups. 
Let $\cal D$ be a collection of finitely generated
groups satisfying the following.

\begin{itemize}
\item Any element $C\in {\cal D}$ has the following properties. 
Quotients and subgroups of $C$
belong to ${\cal D}$. $C$ is residually finite and the FICwF$^P$ is true for
the mapping torus $C\rtimes \l t\r$ for any action of the
infinite cyclic group $\l t\r$ on $C$.
\end{itemize} 

Assume that the vertex groups of $\cal G$ belong to $\cal D$.
Then the FICwF$^P$ is true for $\pi_1({\cal G})$ if ${\cal G}$ satisfies the
intersection property.
\end{thm}

Theorem \ref{ip} is a special case of Proposition \ref{residually-prop} $(2)$ 
and is proved in Section 7.

As a consequence of Theorem \ref{ip} we prove the following.

Let us first recall that a group $G$ is called {\it subgroup separable} 
if the following is satisfied. 
For any finitely generated subgroup $H$ of $G$ and $g\in G-H$  
there is a finite index normal
subgroup 
$K$ of $G$ so that $H\subset K$ and $g\in G-K$. Equivalently, 
a group is subgroup separable if the finitely generated subgroups 
of $G$ are closed in the profinite topology of $G$. 
A subgroup separable group is 
therefore residually finite. The classes of subgroup separable and 
residually finite groups
are extensively studied in group theory and 
Manifold Topology.

\begin{defn} \label{almostdefn} {\rm Let $\cal G$ be a graph 
of groups. An edge $e$ of $\cal G$ is called a {\it finite edge} 
if the edge group ${\cal G}_e$ is finite. $\cal G$ is called 
{\it almost a tree 
of groups} if there are finite edges $e_1,e_2,\ldots$ so that the  
components of ${\cal
G}-\{e_1,e_2,\ldots \}$ are trees. If we remove all the finite edges 
from a graph of groups then the components of the resulting 
graph are called {\it component subgraphs}.}\end{defn}

\begin{thm} \label{introthm} Let $\cal G$ be a graph of finitely
generated groups. Then the FICwF$^P$ is true for $\pi_1({\cal G})$ if one of the
following five conditions is satisfied.

$(1).$ The vertex groups are virtually polycyclic and ${\cal G}$ satisfies
the intersection property.

$(2).$ The vertex groups of $\cal G$ are finitely generated nilpotent and 
$\pi_1({\cal G})$ is subgroup separable.

$(3).$ The vertex and the edge groups of
any component subgraph (see Definition \ref{almostdefn}) are
fundamental
groups of closed surfaces of genus $\geq 2$. Given a component 
subgraph $\cal H$ which has at least one edge there is a subgroup 
$C<\cap_{e\in E_{\cal H}}{\cal H}_e$ 
which is finite index in some edge group and is normal in $\pi_1({\cal H})$.

$(4).$ The vertex groups are virtually cyclic and 
the graph is almost a tree.
Further assume that either the edge groups are finite or the
infinite vertex groups are abelian.

$(5).$ $\cal G$ is almost a tree of groups  
and the 
vertex and the edge groups of any component subgraph of $\cal G$   
are finitely generated abelian and 
of the same rank.\end{thm}

For examples of graphs of groups satisfying the hypothesis 
in $(3)$ see Example \ref{exm}. 

In the next section we start with stating the general 
Fibered Isomorphism conjecture for equivariant 
homology theory (\cite{BL}). We 
show that under certain conditions a group acting on a tree satisfies this 
general conjecture provided the stabilizers also satisfy the 
conjecture. Then we restrict to the pseudoisotopy case of the
conjecture and prove the Theorems. 

We work in this general setting as we will apply these 
methods in the $L$-theory version of the conjecture in later 
works (\cite{R4}, \cite{R5}). We will see that using these methods 
we can prove some stronger results in the $L$-theory case.

Finally, in Examples \ref{end}, 
\ref{end0} and \ref{end1} 
we provide explicit examples of groups for which 
the results of this paper can be applied to 
prove the Fibered Isomorphism conjecture in the 
pseudoisotopy case and show that these groups were not considered 
before. We further show that the 
groups in these examples are neither $CAT(0)$ nor hyperbolic. 
Here we note that in a recent paper 
(\cite{BL1}) Bartels and L\"{u}ck have proved the Fibered 
Isomorphism conjectures in $L$- and lower $K$-theory 
case for hyperbolic groups and for $CAT(0)$-groups which act  
on finite dimensional $CAT(0)$-spaces.

\noindent
{\bf Acknowledgment.}
This work started during my visit to the 
Mathematisches Institut, Universit\"{a}t M\"{u}nster, M\"{u}nster,
Germany under a 
Fellowship of the Alexander von Humboldt-Stiftung. I would like to thank 
my host Wolfgang L\"{u}ck for making this visit possible and  
Arthur Bartels for some useful comments 
during the initial stage of this work. I am also grateful to F. T. Farrell 
for some helpful email communications. Finally thanks to Henrik Rueping 
for pointing out an error in an earlier version of the paper.

\section{Statements of the conjecture and some propositions}

In this section we recall the statement of the Isomorphism 
conjecture for equivariant homology theories (see [\cite{BL}, Section 1]) and 
state some propositions. 

Let ${\cal H}^?_*$ be an equivariant homology theory with values in 
$R$-modules for $R$ a commutative associative ring with unit. 
An equivariant homology theory assigns to 
a group $G$ a $G$-homology theory ${\cal H}^G_*$ which, for a pair 
of $G$-CW  complex $(X,A)$, produces a ${\Bbb Z}$-graded $R$-module 
${\cal H}^G_n(X,A)$. For details see [\cite{L1}, Section 1].

A {\it family} of subgroups of a group $G$ is defined as a set   
of subgroups of $G$ which is closed under taking subgroups and 
conjugations. If $\cal C$ is a class of 
groups which is closed under isomorphisms and taking subgroups then we 
denote by ${\cal C}(G)$ the set of all subgroups of $G$ which belong to 
$\cal C$. 
Then ${\cal C}(G)$ is a family of subgroups of $G$. For example $\cal 
{VC}$, the class of virtually cyclic groups, is closed under isomorphisms 
and taking subgroups. By definition a virtually cyclic group has a 
cyclic subgroup of finite index. Also ${\cal {FIN}}$, the class of 
finite groups is closed under isomorphisms and taking subgroups.
   
Given a group homomorphism $\phi:G\to H$ and a family $\cal C$ of 
subgroups of $H$ define $\phi^*{\cal C}$ to be the family 
of subgroups $\{K<G\ |\ \phi (K)\in {\cal C}\}$ of $G$. Given a family 
$\cal C$ of subgroups of a group $G$ there is a $G$-CW complex $E_{\cal 
C}(G)$ which is unique up to $G$-equivalence satisfying the property that 
for $H\in {\cal C}$ the fixpoint set $E_{\cal C}(G)^H$ is 
contractible and $E_{\cal C}(G)^H=\emptyset$ for $H$ not in ${\cal C}$.

\begin{fjic}([\cite{BL}, Definition 1.1]) Let ${\cal H}^?_*$ be an 
equivariant homology theory with values in $R$-modules. Let $G$ be a group 
and $\cal C$ be a family of subgroups of $G$. Then the {\it 
Isomorphism conjecture} for the pair $(G, {\cal C})$ states that the 
projection 
$p:E_{\cal C}(G)\to pt$ to the point $pt$ induces an isomorphism 
$${\cal H}^G_n(p):{\cal H}^G_n(E_{\cal C}(G))\simeq {\cal H}^G_n(pt)$$ for $n\in 
{\Bbb Z}$. 

And the {\it Fibered Isomorphism conjecture} for the pair $(G, {\cal 
C})$ states that for any group homomorphism $\phi: K\to G$ the 
Isomorphism conjecture is true for the pair $(K, \phi^*{\cal C})$.\end{fjic}

Let $\cal C$ be a class of groups which is closed under isomorphisms and 
taking subgroups.

\begin{defn} \label{definition} {\rm If the 
(Fibered) Isomorphism conjecture is
true for the pair $(G, {\cal C}(G))$ we say that the 
{\it (F)IC$_{\cal C}$ is true for $G$} or simply say 
{\it (F)IC$_{\cal C}(G)$ is satisfied}. Also we say that the {\it 
(F)ICwF$_{\cal C}(G)$ is satisfied} if
the {\it (F)IC$_{\cal C}$} is true for $G\wr H$ for any finite group
$H$. Finally a
group homomorphism $p:G\to K$ is said to {\it satisfy the FIC$_{\cal C}$}  or {\it the 
FICwF$_{\cal C}$} 
if for $H\in p^*{\cal C}(K)$ the FIC$_{\cal C}$ or the FICwF$_{\cal C}$ is true for $H$ respectively.}\end{defn} 

Let us denote by $P$, $K$, $L$, and 
$KH$ the equivariant homology theories arise in  
the stable topological pseudoisotopy theory, 
 algebraic $K$-theory, $L$-theory and in the homotopy $K$-theory 
(\cite{BL}) respectively. The conjectures corresponding to these   
equivariant homology theories are denoted by 
(F)IC$^X_{\cal C}$ (or FICwF$^X_{\cal C}$) where 
$X=P,K,L$ or $KH$ respectively. For the first three theories 
we shorten the notation (F)IC$^X_{\cal {VC}}$ (or FICwF$^X_{\cal {VC}}$) 
to (F)IC$^X$ (or FICwF$^X$). And we shorten (F)IC$^{KH}_{\cal {FIN}}$ 
(or FICwF$^{KH}_{\cal {FIN}}$) to (F)IC$^{KH}$ (or FICwF$^{KH}$). 
The Isomorphism conjectures (F)IC$^P$, (F)IC$^K$ and 
(F)IC$^L$ are equivalent to the Farrell-Jones conjectures 
stated in (\S 1.7) \S 1.6 in \cite{FJ}. (For details see [\cite{BL}, Sections 5 and 6] for
the $K$ and $L$ theories and see [\cite{LR}, 4.2.1 and 4.2.2] for the pseudoisotopy theory.) 

\begin{defn} \label{property} {\rm We say that {\it ${\cal T}_{\cal C}$} 
({\it $_w{\cal T}_{\cal C}$}) is 
satisfied if for a
graph of groups $\cal G$ with vertex groups 
from the class $\cal C$ 
the FIC$_{\cal
C}$ (FICwF$_{\cal
C}$) for $\pi_1({\cal G})$ is true. 

Let us now assume that $\cal C$ contains all the 
finite groups. We say that {\it
$_t{\cal T}_{\cal C}$} 
({\it $_f{\cal T}_{\cal C}$})
is satisfied if for a graph of 
groups $\cal G$ with trivial (finite) edge
groups and vertex
groups belonging to the class $\cal C$, the FIC$_{\cal
C}$ for $\pi_1({\cal G})$ is true. If we replace the FIC$_{\cal
C}$ by the FICwF$_{\cal C}$  
then we denote the corresponding properties by {\it
$_{wt}{\cal T} _{\cal C}$} (
{\it $_{wf}{\cal T} _{\cal C}$}). Clearly 
$_{w}{\cal T} _{\cal C}$ implies 
${\cal T} _{\cal C}$ and $_{w*}{\cal T} _{\cal C}$ implies 
$_{*}{\cal T} _{\cal C}$ where $*=t$ or 
$f$.

And we say that {\it ${\cal P} _{\cal C}$} is 
satisfied if for $G_1, G_2\in {\cal C}$ the product
$G_1\times G_2$ satisfies the FIC$ _{\cal C}$.

{\it ${\cal {FP}} _{\cal C}$} is satisfied if 
whenever the FIC$ _{\cal C}$ is true for two groups 
$G_1$ and $G_2$ then the FIC$ _{\cal C}$ is true for 
the free product $G_1*G_2$.

We denote the above properties for the  
equivariant homology theories $P,K, L$ or $KH$  with only a 
super-script by $P,K, L$ or $KH$ respectively. For example, 
${\cal T}_{\cal C}$ 
for $P$ is denoted by ${\cal T}^P$ since in all the first three 
cases we set ${\cal C}={\cal {VC}}$ and for $KH$  we 
set ${\cal C}={\cal {FIN}}$.}\end{defn}

We will show in this article and in \cite{R4} 
that the above properties are satisfied in several 
instances of 
the conjecture. 

For the rest of this section we assume that $\cal C$ is also closed under 
quotients and contains all the finite groups. 

\begin{prop} \label{theorem} {\bf (Graphs of groups).} 
Let ${\cal C}={\cal {FIN}}$ or $\cal {VC}$. Let $\cal G$ be a graph of
groups and there is a homomorphism 
$f:\pi_1({\cal G})\to Q$ onto another group $Q$ so that the restriction
of $f$ to any vertex group has finite kernel. 
If the FIC$ _{\cal C}(Q)$ 
(FICwF$ _{\cal C}(Q)$) is satisfied then 
the FIC$ _{\cal C}(\pi_1({\cal G}))$ 
(FICwF$ _{\cal C}(\pi_1({\cal G}))$) is 
also satisfied
provided one of the 
following holds.
(In the FICwF$ _{\cal C}$-case in addition assume that 
${\cal P} _{\cal C}$ is 
satisfied.)

$(a).$ ${\cal T} _{\cal C}$ 
($_w{\cal T} _{\cal C}$) is satisfied.

$(b).$ The edge groups of $\cal G$ are finite and 
$_f{\cal T} _{\cal C}$ 
($_{wf}{\cal T} _{\cal C}$) is satisfied.

$(c).$ The edge groups of $\cal G$ are trivial and 
$_t{\cal T} _{\cal C}$ 
($_{wt}{\cal T} _{\cal C}$) is satisfied.

\end{prop}

Proposition \ref{theorem} will be proved in Section 5.

For definition of continuous ${\cal H}^?_*$ in the following 
statement 
see [\cite{BL}, Definition 3.1]. Also see Proposition \ref{proposition}.

\begin{prop} \label{residually-prop}{\bf (Graphs of residually finite
groups)} Assume that ${\cal P} _{\cal C}$ and 
$_{wt}{\cal T} _{\cal C}$ are satisfied. Let $\cal
G$ be a graph of groups. If $\cal G$ is infinite then assume that 
${\cal H}^?_*$ is continuous. 

$(1).$ Assume that the edge groups of $\cal G$ are finite and the vertex
groups are residually finite. If the FICwF$ _{\cal {C}}$ is
true for the vertex groups of $\cal G$, then it is true for $\pi_1({\cal G})$.

For the next three items assume that ${\cal C}={\cal {VC}}$.

$(2).$ Let $\cal D$ be the collection of groups as defined in Theorem
\ref{ip} replacing `FICwF$^P$' by `FICwF$ _{\cal
{VC}}$'. Assume that the vertex groups of $\cal G$ belong to $\cal D$. Then the 
FICwF$ _{\cal {VC}}$ is true for $\pi_1({\cal G})$ if $\cal
G$ satisfies the intersection property.

$(3).$ Assume that the vertex groups of $\cal G$ are virtually polycyclic
and the FICwF$ _{\cal {VC}}$ is true for virtually
polycyclic groups. Then the FICwF$ _{\cal {VC}}$ is true for
$\pi_1({\cal G})$ provided either $\cal G$ satisfies the
intersection property or the vertex groups are finitely generated
nilpotent and $\pi_1({\cal G})$ is subgroup separable.

$(4).$ Assume that the vertex 
and the edge groups of
any component subgraph (see Definition \ref{almostdefn}) are
fundamental
groups of closed surfaces of genus $\geq 2$ and for every component 
subgraph $\cal H$ which has at least one edge there is a subgroup 
$C<\cap_{e\in {\cal H}}{\cal H}_e$ 
which is of finite index in some edge group and is normal in $\pi_1({\cal H})$.
Then 
the FICwF$ _{\cal {VC}}$ is true for $\pi_1({\cal G})$
provided the FICwF$ _{\cal {VC}}$ is true for the
fundamental groups of closed $3$-manifolds which fibers over the circle.

\end{prop}

Let $F^{\infty}$ denote a countable infinitely generated free group and 
${\Bbb Z}^{\infty}$ denotes a countable infinitely generated abelian 
group. Also 
$G\rtimes H$ denotes a semidirect product with 
respect to an arbitrary action of $H$ on $G$. When 
$H$ is infinite cyclic and generated by the symbol $t$, 
we denote it by $\l t\r$

\begin{prop} \label{maintheorem1} {\bf (Graphs of abelian groups).} 
Let $\cal G$ be a graph of groups whose vertex groups are 
finitely generated abelian and let ${\cal H}^?_*$ be continuous. 

$(1).$ Assume that the FIC$ _{\cal
{VC}}$ is true for $F^{\infty}\rtimes \l
t\r$ and that ${\cal P} _{\cal
{VC}}$ is satisfied.
Then the FIC$_{\cal 
{VC}}$ is true for $\pi_1({\cal G})$ provided one of the 
following holds.

(a). $\cal G$ is a tree and the vertex groups of 
$\cal G$ are torsion free.  

For the next two items assume that the FICwF$ _{\cal
{VC}}$ is true for $F^{\infty}\rtimes \l
t\r$.

(b). $\cal G$ is a tree.

(c). $\cal G$ is not a tree and 
the FIC$ _{\cal 
{VC}}$ is true for ${\Bbb Z}^{\infty}\rtimes \l t\r$ for any countable 
infinitely generated abelian group ${\Bbb Z}^{\infty}$. 

$(2).$ Assume that ${\cal P} _{\cal {VC}}$ and 
$_{wt}{\cal T} _{\cal {VC}}$ are satisfied. 
Further assume that $\cal G$ is almost a tree of groups and the
vertex and the edge groups of any component subgraph of $\cal G$ 
have the same rank. Then the FICwF$_{\cal
{VC}}$ is true for $\pi_1({\cal G})$ provided one of the followings 
is satisfied.

(i). The FICwF$ _{\cal
{VC}}$ is true for ${\Bbb Z}^{\infty}\rtimes \l t\r$ for any countable
infinitely generated abelian group ${\Bbb Z}^{\infty}$.

(ii). The vertex and the edge groups of any component subgraph of 
$\cal G$ have rank equal to $1$.

$(3).$ The IC$^{L}_{\cal {FIN}}$ is true for 
$\pi_1({\cal 
G})$ provided the vertex groups of $\cal G$ are torsion free.
\end{prop}

In $(1)$ of Proposition \ref{maintheorem1}
 if we assume that the  FICwF$ _{\cal 
{VC}}$ is true for ${\Bbb Z}^{\infty}\rtimes \l t\r$ and for 
$F^{\infty}\rtimes \l t\r$ then one can deduce from the same 
proof, using $(3)$ of Proposition \ref{product} instead 
of Lemma \ref{inverse}, that the  FICwF$ _{\cal
{VC}}$ is true for $\pi_1({\cal G})$ irrespective of whether 
$\cal G$ is a tree or not.

\begin{prop} \label{TC} {\bf ($_w{\cal T}^{P}$).} Let $\cal G$ be a 
graph of virtually cyclic groups so that the graph is almost a tree. 
Further assume that either the edge groups are finite or the 
infinite vertex groups are abelian.
Then the FICwF$^{P}$ is true for 
$\pi_1({\cal G})$.\end{prop}

An immediate corollary of the Proposition is the following.

\begin{cor} \label{tt} {\bf ($_{wf}{\cal T}^{P}$).} $_{wf}{\cal T}^{P}$ 
(and hence $_{wt}{\cal T}^{P}$) is satisfied.\end{cor}

\begin{rem}\label{possible}{\rm We remark here that it 
is not yet known if the FIC$^{P}$ is true 
for the HNN-extension $G=C*_C$, with 
respect to the maps $id:C\to C$ and $f:C\to C$ where $C$ is an infinite 
cyclic group, $id$ is the identity map and $f(u)=u^2$ for $u\in C$. This 
was already mentioned in the introduction in \cite{FL}. Note 
that $G$ is isomorphic to the semidirect product ${\Bbb 
Z}[\frac{1}{2}]\rtimes \l t\r$, where $t$ acts on ${\Bbb 
Z}[\frac{1}{2}]$ by multiplication by $2$. The main problem with this
example is that ${\Bbb 
Z}[\frac{1}{2}]\rtimes \l t\r$ is not subgroup separable.}\end{rem}

Now we recall that a property ({\it tree property}) similar to ${\cal 
T} _{\cal C}$ was defined in 
[\cite{BL}, Definition 4.1]. The 
tree property of \cite{BL} is stronger than ${\cal T}_{\cal 
C}$. Corollary 4.4 in \cite{BL} was proved under the assumption 
that the tree property is satisfied. In the following proposition we state  
that this is true with the weaker assumption that ${\cal 
T} _{\cal 
{FIN}}$ is satisfied. The proof of the Proposition goes exactly in 
the same way as the proof of [\cite{BL}, Corollary 4.4]. It can also be deduced
from Proposition \ref{theorem}.

\begin{prop} \label{finite} Let 
$1\to K\to G\to Q\to 1$ be an exact sequence of groups. Assume that  
${\cal T} _{\cal {FIN}}$ is satisfied and $K$ acts 
on a tree with 
finite stabilizers and that the FIC$ _{\cal {FIN}}(Q)$ is 
satisfied. Then the FIC$ _{\cal {FIN}}(G)$ is also 
satisfied.\end{prop} 

\section{Graphs of groups}

In this section we prove some results on graphs of groups needed for the 
proofs of the Theorems and Propositions.

We start by recalling that by Bass-Serre theory 
a group acts on a tree without inversion if and only 
if the group is isomorphic to the 
fundamental group of a graph of 
groups (the Structure Theorem I.4.1 in \cite{DD}). Therefore,  
throughout the paper by `action of a group on a tree' we will mean 
an action without inversion.

\begin{lemma} \label{almost} Let $\cal G$ be a finite and almost a tree
of groups (see Definition \ref{almostdefn}). 
Then there is another graph of groups $\cal H$ with the following 
properties.

$(1).$ $\pi_1({\cal G})\simeq \pi_1({\cal H})$.

$(2).$ Either $\cal H$ has no edge or the edge groups 
of $\cal H$ are finite.

$(3).$ The vertex groups of $\cal H$ are of the form $\pi_1({\cal K})$ 
where $\cal K$ varies over subgraphs of groups of $\cal G$ which are 
maximal with respect to the 
property that the underlying graph of $\cal K$ is a tree and the 
edge (if there is any) groups of 
$\cal K$ are all infinite.\end{lemma}

\begin{proof} 
The proof is by
induction on the number of finite edges of the graph. Recall 
that an edge is called a finite edge if the corresponding edge 
group is finite (Definition \ref{almostdefn}).
If the graph 
has no finite edge then by definition 
of almost a graph of groups $\cal G$ is a tree. In this case we can take 
$\cal H$ to be a graph consisting of a single vertex and associate the 
group $\pi_1({\cal G})$ to the vertex. So assume 
${\cal G}$ has
$n$ finite edges and that the Lemma is true for graphs with $\leq 
n-1$ finite edges.   
Let $e$ be an edge of ${\cal G}$ with ${\cal G}_e$ finite. 
If ${\cal G}-\{e\}$ is connected
then $\pi_1({\cal G})\simeq \pi_1({\cal
G}_1)*_{{\cal G}_e}$ where ${\cal G}_1={\cal G}-\{e\}$ is a graph 
with $n-1$ finite 
edges. By the induction hypothesis there 
is a graph of groups ${\cal H}_1$ satisfying $(1), (2)$ and $(3)$ for 
${\cal G}_1$. Let $v_1$ and $v_2$ be the vertices of $\cal G$ to 
which the ends of $e$ are 
attached and let $v'_1$ and $v'_2$ be the vertices of 
${\cal H}_1$ so that ${\cal G}_{v_i}$ is a 
subgroup of ${\cal H}_{v'_i}$ for $i=1,2$. (Note 
that $v_1$ and $v_2$ could be the same vertex). 
Define $\cal H$ by attaching an edge 
$e'$ to ${\cal H}_1$ so that the ends of $e'$ are attached to 
$v'_1$ and $v'_2$ and associate 
the group ${\cal G}_e$ to $e'$. The injective homomorphisms  
${\cal H}_{e'}\to {\cal H}_{v'_i}$ for $i=1,2$   
are defined by the homomorphisms ${\cal G}_e\to {\cal G}_{v_i}$.  It is now 
easy to check that $\cal H$ satisfies $(1), (2)$ and $(3)$ for $\cal G$.
On the other hand if ${\cal G}-\{e\}$ has two components say ${\cal G}_1$ and 
${\cal G}_2$ then $\pi_1({\cal G})\simeq \pi_1({\cal
G}_1)*_{{\cal G}_e} \pi_1({\cal G}_2)$ where ${\cal G}_1$ and ${\cal G}_2$ 
has $\leq
n-1$ finite edges. Using the induction hypothesis and 
a similar argument as above we complete the proof of the  
Lemma.\end{proof}

\begin{lemma}\label{finitefree} A finitely generated group contains a free 
subgroup of finite 
index if and only if the group acts on a tree with finite 
stabilizers. And a group acts on a tree with trivial 
stabilizers if and only if the group is free.\end{lemma}

\begin{proof} By  
[\cite{DD}, Theorem IV.1.6] a group contains a free subgroup of finite index if 
and only if the group acts on a tree with finite stabilizers 
and the stabilizers have bounded order. The proof of the Lemma 
now follows easily.\end{proof}

We will also need the following two Lemmas.

\begin{lemma} \label{stark} Let $\cal G$ be a graph of finitely 
generated abelian groups so that the underlying graph of $\cal G$ is a 
tree. Then the restriction of the abelianization 
homomorphism $\pi_1({\cal G})\to H_1(\pi_1({\cal G}), {\Bbb Z})$ to each 
vertex group of 
the tree of groups $\cal G$ is injective.\end{lemma}

\begin{proof} If the tree $\cal G$ is finite and the vertex groups are 
finitely generated free abelian then it was proved in 
[\cite{S}, Lemma 3.1] that there is a homomorphism $\pi_1({\cal 
G})\to A$ onto a free abelian group so that the restriction of this 
homomorphism to each vertex group is injective. In fact it was shown  
there that the abelianization homomorphism $\pi_1({\cal 
G})\to H_1(\pi_1({\cal G}), {\Bbb Z})$ is injective when restricted 
to the vertex groups  
and then since the vertex groups are torsion free 
$\pi_1({\cal
G})\to A=H_1(\pi_1({\cal 
G}), {\Bbb Z})/\{torsion\}$ is also injective on the vertex groups. 
The same proof 
goes 
through, without the torsion free vertex group assumption, to 
prove the Lemma when $\cal G$ is finite. Therefore we 
mention 
the additional arguments needed in the infinite case. At first write 
${\cal G}$ as 
an increasing union of finite trees ${\cal G}_i$ of finitely 
generated abelian groups. Consider the following commutative diagram.

$$\begin{CD}
\pi_1({\cal G}_i) @>>> H_1(\pi_1({\cal 
G}_i), {\Bbb Z})\\
@VVV    @VVV\\
\pi_1({\cal G}_{i+1}) @>>> H_1(\pi_1({\cal 
G}_{i+1}), {\Bbb Z})
\end{CD}$$

Note that the left hand side vertical map is injective. And by the finite 
tree case the restriction of the two horizontal maps to each vertex 
group of the respective trees are injective.
Now since group homology and fundamental group commute 
with direct limit, taking limit completes the proof 
of the Lemma.\end{proof}

\begin{lemma} \label{fi} Let $\cal G$ be a finite graph of finitely 
generated groups
satisfying the following.

{\bf P.} Each edge group 
is of finite index in the end vertex groups of the
edge. Also assume that the intersection of the edge 
groups contains a subgroup $C$ (say) which is normal 
in $\pi_1({\cal G})$ and is of finite index in some edge ($e$ say) group. 

Then $\pi_1({\cal
G})/C$ is isomorphic to the fundamental group of a  
finite-type (see Definition \ref{def1.1}) 
graph of groups. Consequently $\cal G$ has the
intersection property.\end{lemma} 

\begin{proof} 
The proof is by induction on the number of edges of the graph. 

{\bf Induction hypothesis.} $(IH_n)$ For any finite graph of groups $\cal 
G$ with 
$\leq n$ edges which satisfies {\bf P}, $\pi_1({\cal
G})/C$ is isomorphic to the fundamental group of a graph of groups whose
underlying graph is the same as that of $\cal G$ and the vertex 
groups are finite and isomorphic to ${\cal G}_v/C$ where 
$v\in V_{\cal G}$. 

If $\cal G$ has one edge $e$ then $C<{\cal G}_e$ and $C$ 
is normal in $\pi_1({\cal G})$. First assume $e$ disconnects 
$\cal G$. Since ${\cal G}_e$ is of 
finite index in the end vertex groups and $C$ is also 
of finite index in ${\cal G}_e$, $C$ is of finite index in the 
end vertex groups of $e$. Therefore $\pi_1({\cal G})/C$ 
has the desired property. The argument is same when 
$e$ does not disconnect $\cal G$.
 
Now assume $\cal G$ has $n$ edges and satisfies {\bf P}.

Let us first consider the case that there is an 
edge $e'$ other than $e$ so that  
${\cal G}-\{e'\}={\cal D}$ (say)  
is a connected graph. Note that $\cal D$ has $n-1$ edges 
and satisfies $\bf P$.

Hence by $IH_{n-1}$, $\pi_1({\cal D})/C$ is isomorphic 
to the fundamental group of a finite-type 
graph of groups with $\cal D$ as the underlying 
graph and the vertex groups 
are of the form ${\cal D}_v/C$.  

Let $v$ be an end vertex of $e'$. Then 
${\cal D}_v/C$ is finite. 
Also by hypothesis ${\cal G}_{e'}$ 
is of finite index in ${\cal G}_v={\cal D}_v$. Therefore 
${\cal G}_{e'}/C$ is also finite. This 
completes the proof in this case.

Now if every edge $e'\neq e$ disconnects $\cal G$ then $\cal G-\{e\}$ is a 
tree. Let $e'$ be an edge other than $e$ so that one 
end vertex $v'$ (say) of $e'$ has valency $1$. Such an edge 
exists because $\cal G-\{e\}$ is a tree. Let 
${\cal D}={\cal G}-(\{e'\}\cup \{v'\})$. 
Then 
$\pi_1({\cal G})\simeq \pi_1({\cal D})*_{{\cal G}_{e'}}{\cal G}_{v'}$. 
Let $v''$ be the other end vertex of $e'$. Then by the induction 
hypothesis $C$ is of finite index in 
${\cal G}_{v''}={\cal D}_{v''}$. Hence $C$ 
is of finite index in ${\cal G}_{e'}$. Also ${\cal G}_{e'}$ is of 
finite index in ${\cal G}_{v'}$. Therefore $C$ is also of finite index 
in ${\cal G}_{v'}$. 

This completes the proof.  
\end{proof}

The following Lemma and Example give some concrete examples of graphs of groups
with intersection property.

\begin{lemma} \label{ic} Let $\cal G$ be a finite graph of groups so that 
all the vertex and the edge groups are finitely generated abelian and of the 
same rank $r$ (say) and the underlying graph of $\cal G$ is a tree. Then 
$\cap_{e\in E_{\cal G}}{\cal G}_e$ contains a rank $r$ free 
abelian subgroup 
$C$ which is normal in $\pi_1({\cal G})$ so that $\pi_1({\cal
G})/C$ is isomorphic to the fundamental group of a graph of groups whose 
underlying graph is $\cal G$ and vertex groups are finite and isomorphic 
to ${\cal G}_v/C$ where $v\in V_{\cal G}$.\end{lemma}

\begin{proof} The proof is by induction on the number of edges. 
If the graph has one edge $e$ then clearly ${\cal G}_e$ is normal 
in $\pi_1({\cal G})$, for ${\cal G}_e$ is normal in the two end 
vertex groups. This follows from Lemma \ref{gene}. 
So by induction assume that the Lemma is true for 
graphs with $\leq n-1$ edges. Let $\cal G$ be a finite graph with $n$ 
edges and satisfies the hypothesis of the Lemma. Consider an 
edge $e$ which has one end vertex $v$ (say) with valency $1$. Such 
an edge exists because the graph is a tree. Let 
$v_1$ be the other end vertex of $e$. Then 
$\pi_1({\cal G})\simeq \pi_1({\cal G}')*_{{\cal G}_e}{\cal G}_v$. 
Here ${\cal G}'={\cal G}-(\{e\}\cup\{v\})$. Clearly by induction hypothesis 
there is a finitely generated free abelian normal subgroup 
$C_1<\cap_{e'\in E_{{\cal G}'}}{{\cal G}'}_{e'}$ of rank 
$r$ of $\pi_1({\cal G}_1)$ satisfying all the required properties. 
Now note that $C_1\cap {\cal G}_e=C'$ (say) is of finite index 
in ${{\cal G}'}_{v_1}$ and also in $C_1$ and $C'$ has rank $r$. 
This follows from the following easy to 
verify Lemma. Now since $C_1$ is finitely generated and $C'$ is of finite 
index in $C_1$ we can find a characteristic subgroup $C (<C')$ of $C_1$ of finite 
index. Therefore $C$ has rank $r$ and is normal in $\pi_1({\cal G}')$, since 
$C_1$ is normal in $\pi_1({\cal G}')$. Obviously $C$ is normal in 
${\cal G}_v$. Therefore, we again use Lemma \ref{gene} 
to conclude that $C$ is normal in 
$\pi_1({\cal G})$. The other properties are clearly satisfied.
This completes the proof of the Lemma.

\begin{lemma} \label{rank} Let $G$ be a finitely generated abelian group 
of rank $r$. Let $G_1,\ldots , G_k$ be rank $r$ subgroups of $G$. Then 
$\cap_{i=1}^{i=k}G_i$ is of rank $r$ and of finite index in 
$G$.\end{lemma}  
\end{proof}

\begin{exm}\label{exm} {\rm Let $G$ and $H$ be finitely generated groups 
and let $H$ be a finite index normal subgroup of $G$. Let $f:H\to G$ 
be the inclusion. Consider a finite tree of groups $\cal G$ whose vertex groups 
are copies of $G$ and the edge groups are copies of $H$. Also assume 
that the maps from the edge groups to the vertex groups (defining the 
tree of group structure) are $f$. Then $\cal G$ has the intersection 
property.}\end{exm}

\begin{lemma} \label{gene} Let $G=G_1*_HG_2$ be a generalized 
free product. If $H$ is normal 
in both $G_1$ and $G_2$ then $H$ is normal in $G$.\end{lemma}

\begin{proof} The proof follows by using the normal form of 
elements in a generalized free product. See [\cite{LS}, p. 72].
\end{proof}

\section{Residually finite groups} 
In this section we recall and also prove some basic results we need on
residually finite groups. For this section we 
abbreviate `residually finite' by $\cal {RF}$.

\begin{lemma} \label{graph-res} The
fundamental group of a finite graph of $\cal {RF}$ groups with finite edge
groups is $\cal {RF}$.\end{lemma}

\begin{proof} The proof is by induction on the number of edges. If there
is no edge then there is nothing prove. So assume the Lemma for graphs
with $\leq n-1$ edges. Let $\cal G$ be a graph of groups with $n$ edges
and satisfies the hypothesis. It follows that $\pi_1({\cal G})\simeq
\pi_1({\cal G}_1)*_F\pi_1({\cal G}_2)$ or $\pi_1({\cal G})\simeq
\pi_1({\cal G}_1)*_F$ where $F$ is a finite group and ${\cal G}_i$ satisfy
the hypothesis of the Lemma and has $\leq n-1$ edges. Also note that
$\pi_1({\cal G}_1)*_F\pi_1({\cal G}_2)$ can be embedded as a 
subgroup in $(\pi_1({\cal
G}_1)*\pi_1({\cal G}_2))*_F$. Therefore by the induction hypothesis and
since
free product of $\cal {RF}$ groups is again
$\cal {RF}$ (see \cite{Gr}) we only need to prove that, for
finite $H$, $G*_H$
is $\cal {RF}$ if so is $G$. But, this follows from
\cite{BT} or  [\cite{DC}, Theorem 2]. This completes the proof 
of the Lemma.\end{proof}

\begin{lemma} \label{res-ext} Let $1\to K\to G\to H\to 1$ be an extension
of groups so that $K$ and $H$ are ${\cal {RF}}$. Assume that any finite
index
subgroup of $K$ contains a subgroup $K'$ so that $K'$ is
normal in $G$ and $G/K'$ is $\cal {RF}$. 
Then $G$ is $\cal {RF}$.\end{lemma}

\begin{proof} Let $g\in G-K$ and $g'\in H$ be the image of $g$ in $H$.
Since $H$ is ${\cal {RF}}$ there is a finite index subgroup $H'$ of $H$
not
containing $g'$. The inverse image of $H'$ in $G$ is a finite index
subgroup not containing $g$. 

Next let $g\in K-\{1\}$. Choose a finite index subgroup of $K$ not containing
$g$. By hypothesis there is a finite index subgroup $K'$ of $K$ which is
normal in $G$ and does not contain $g$. Also $G/K'$ is ${\cal {RF}}$. Now
applying the previous case we complete the proof.\end{proof}

\begin{lemma} \label{inter} Let $\cal G$ be a finite graph of 
finitely generated $\cal {RF}$  
groups satisfying the intersection property. Assume the following. 
`Given an edge $e$ and an end vertex $v$ of $e$, 
for every subgroup $E$ of ${\cal G}_e$ 
which is normal in ${\cal G}_v$, the quotient ${\cal G}_v/E$ is again $\cal {RF}$.' 
Then $\pi_1({\cal
G})$ is $\cal {RF}$.\end{lemma}

\begin{proof} Using Lemma \ref{graph-res} we can assume that all the edge
groups of $\cal G$ are infinite. Now the proof is by induction on the number
of edges of the graph $\cal G$. Clearly the induction starts because 
if there is no edge then the Lemma is true.
Assume the result for all graphs satisfying the hypothesis with
number of edges $\leq n-1$ and let $\cal G$ be a graph of groups with $n$
number of edges and satisfying the hypothesis of the Lemma. By the 
intersection property 
there is a normal subgroup $K$, contained 
in all the edge groups, of $\pi_1({\cal G})$ which is of finite index in
some edge group. Hence we have the following. 

{\bf K.} 
$\pi_1({\cal G})/K$ is isomorphic to the fundamental group of a finite
graph of $\cal {RF}$ groups (by hypothesis the quotient of a vertex group 
by $K$ is $\cal {RF}$) with $n$ edges and some edge 
group is finite. Also it is easily seen that this 
quotient graph of groups has the intersection property. 

By the induction hypothesis and by Lemma \ref{graph-res} $\pi_1({\cal
G})/K$ is $\cal {RF}$. 

Now we would like to apply Lemma \ref{res-ext} to the exact sequence.

$$1\to K\to \pi_1({\cal G})\to \pi_1({\cal G})/K\to 1.$$

Let $H$ be a finite index subgroup of $K$. Since $K$ is finitely
generated (being of finite index in a finitely generated group) we can
find a finite index characteristic subgroup $H'$ of $K$ contained in $H$.
Hence $H'$ is normal in $\pi_1({\cal G})$. It
 is now easy to see that {\bf K} is satisfied if we replace $K$ by $H'$.
Hence $\pi_1({\cal G})/H'$ is $\cal {RF}$.

Therefore, by Lemma \ref{res-ext} $\pi_1({\cal G})$ is $\cal {RF}$.
\end{proof}

\begin{lemma} \label{rfcyclic} Let $\cal G$ be a finite graph of
virtually cyclic groups
so that either the edge groups are finite or the 
infinite vertex groups are abelian and 
the associated graph is almost a tree. Then $\pi_1({\cal 
G})$ is ${\cal {RF}}$.\end{lemma}

\begin{proof} Applying Lemma \ref{graph-res} and using the 
definition of almost a graph of groups we can assume 
that the graph is a tree and all the edge groups are infinite. 
Now using Lemma \ref{ic} we see that the hypothesis of Lemma \ref{inter} is
satisfied. This proves the Lemma.\end{proof}  

\begin{lemma} \label{res-surface} Let $\cal G$ be a finite graph of groups
whose 
vertex and edge groups are fundamental groups of closed 
surfaces of genus $\geq 2$. Also 
assume that the intersection of the edge groups contains a subgroup 
$C$ (say) which 
is normal in $\pi_1({\cal G})$ and is of finite index in some edge group. Then 
$\pi_1({\cal G})$ contains a normal subgroup isomorphic to 
the fundamental group of a closed surface so that the
quotient is isomorphic to the fundamental group of a finite-type  
graph of 
groups. Also
$\pi_1({\cal G})$ is $\cal {RF}$.\end{lemma}

\begin{proof} We need the
following Lemma.

\begin{lemma} Let $S$ be a closed surface.
Let $G$ be a subgroup of $\pi_1(S)$. Then $G$ is isomorphic to the
fundamental group of a closed surface if and only if $G$ is of 
finite index in $\pi_1(S)$.\end{lemma}

\begin{proof} The proof follows from covering space theory.\end{proof}

Therefore, using the above Lemma we get that the edge groups 
of $\cal G$ are of finite index in the
end vertex groups of the corresponding edges. Hence by 
Lemma \ref{fi} $\pi_1({\cal G})/C$
is the fundamental group of a finite-type graph of groups.
This proves the first statement. 

Now 
by Lemma \ref{graph-res} $\pi_1({\cal G})/C$ is $\cal {RF}$. Next, by
\cite{GB1} closed surface groups are $\cal {RF}$. Then using the 
Lemma above, it is
easy to check that the hypothesis of Lemma \ref{res-ext} is satisfied for
the following exact sequence. $$1\to C\to \pi_1({\cal G})\to \pi_1({\cal
G})/C\to 1.$$

Hence $\pi_1({\cal G})$ is $\cal {RF}$.
\end{proof}

\section{Basic results on the Isomorphism conjecture}

In this section we recall some known facts as well as deduce some 
basic results on the Isomorphism conjecture. 

$\cal C$ always denotes a class of groups closed under isomorphisms and 
taking subgroups unless otherwise mentioned.

We start by noting that if the FIC$ _{\cal C}$ is 
true for a group $G$ then the FIC$ _{\cal C}$ is also true 
for any subgroup $H$ 
of $G$. We will refer to this fact as the {\it hereditary property} in this 
paper.

By the Algebraic Lemma in \cite{FR} if $G$ is a 
normal subgroup of $K$ then $K$ can be embedded in 
the wreath product $G\wr (K/G)$. We will be using this fact throughout the 
paper without explicitly mentioning it. 

\begin{lemma} \label{finiteindex} If the FICwF$ _{\cal 
C}(G)$ is satisfied 
then the FICwF$ _{\cal C}(L)$ is also satisfied for any 
subgroup $L$ of $G$. 
\end{lemma} 

\begin{proof} Note that given a group $H$, $L\wr H$ is a subgroup of $G\wr H$. Now use the 
hereditary property of the FIC$ _{\cal
C}$.
\end{proof}

\begin{prop} \label{proposition} Assume that 
${\cal H}^?_*$ is continuous. Let $G$ be a group and 
$G=\cup_{i\in I}G_i$ where $G_i$'s are increasing sequence of subgroups of 
$G$ so that the FIC$ _{\cal C}(G_i)$ is satisfied for  
$i\in I$. Then the FIC$ _{\cal C}(G)$ is also 
satisfied. And if the FICwF$ _{\cal C}(G_i)$ is 
satisfied for 
$i\in I$ then the FICwF$ _{\cal C}(G)$ is also
satisfied.\end{prop}

\begin{proof} The first assertion is the same as the conclusion of 
[\cite{BL}, Proposition 3.4]  
and the second one is easily deducible from it, since given a group $H$, 
$G\wr H=\cup_{i\in I}(G_i\wr H)$.\end{proof}

\begin{rem} \label{finitegraph} {\rm 
Since the fundamental group of an infinite graph of groups can be written 
as an increasing union of fundamental groups of finite subgraphs, 
throughout rest of the paper we consider only finite graphs. The infinite case 
will then follow if the corresponding equivariant homology theory satisfies  
the assumption of Proposition \ref{proposition}. Examples of such 
equivariant homology theories are $P$, $K$, 
$L$, and $KH$. See Section 8.}\end{rem} 

The following Lemma from \cite{BL} is crucial for proofs of the
results in this paper. This result in the context of the 
original Fibered Isomorphism conjecture (\cite{FJ}) 
was proved in [\cite{FJ}, Proposition 2.2].

\begin{lemma}\label{inverse} ([\cite{BL}, Lemma 2.5]) Let $\cal C$ be also 
closed under taking quotients. Let $p:G\to Q$ 
be a surjective group homomorphism and assume that 
the FIC$_{\cal C}$ is true for $Q$ and for $p$. 
Then $G$ 
satisfies the FIC$ _{\cal C}$.\end{lemma} 

\begin{prop} \label{product} Let $\cal C$ be as in the 
statement of the above Lemma. Assume that ${\cal P}_{\cal 
C}$ is 
satisfied. 

(1). If the FIC$ _{\cal C}$ (FICwF$_{\cal C}$) is true 
for $G_1$ and $G_2$ then the FIC$ _{\cal C}$ (
FICwF$_{\cal C}$) is true for $G_1\times G_2$. 

(2). Let $G$ be a finite index normal subgroup of a group $K$. If the  
FICwF$_{\cal C}(G)$ is satisfied then the FICwF$ _{\cal 
C}(K)$ is also satisfied.

(3). Let $p:G\to Q$ be a group homomorphism. If 
the FICwF$_{\cal C}$ is true for $Q$ and for $p$ 
then the FICwF$_{\cal C}$ is true for $G$. 
\end{prop}

\begin{proof} The proof of $(1)$ is essentially two applications of Lemma 
\ref{inverse}. First apply it to the projection $G_1\times G_2\to G_2$. 
Hence to prove the Lemma we need to check that the 
FIC$_{\cal C}$ is true for $G_1\times H$ for 
any $H\in {\cal C}(G_2)$. Now fix $H\in {\cal C}(G_2)$ and apply Lemma 
\ref{inverse} to the projection $G_1\times H\to G_1$. Thus we need to 
show that the FIC$_{\cal C}$ is true for 
$K\times H$ where $K\in {\cal C}(G_1)$. But this is 
exactly ${\cal P} _{\cal C}$ which is true by 
hypothesis. 

Next note that given a group $H$, $(G_1\times G_2)\wr H$ is a 
subgroup of $(G_1\wr H)\times (G_2\wr H)$. Therefore the 
FICwF$_{\cal C}$ is true for $G_1\times G_2$ if it is true for $G_1$ and $G_2$.

For $(2)$ let $H=K/G$. Then $K$ is a 
subgroup of $G\wr H$. Let $L$ be a finite group then it is easy to check 
that $$K\wr L < (G\wr H)\wr L \simeq G^{H\times L}\rtimes (H\wr L) $$
$$< G^{H\times L}\wr (H\wr L) < \Pi_{{|H\times 
L|}-times}(G\wr (H\wr L)).$$

The isomorphism in the above display follows from [\cite{FL}, Lemma 2.5] with respect to 
some action of $H\wr L$ on $G^{H\times L}$ (replace $X$ by $A$ and 
$Y$ by $B$ in [\cite{FL}, Lemma 2.5]). The second inclusion follows from 
the Algebraic Lemma in \cite{FR}.

Now using $(1)$ and by hypothesis we complete the proof.

For $(3)$ we need to prove that the 
FIC$_{\cal C}$ is true for $G\wr H$ for any finite group $H$. We now 
apply Lemma \ref{inverse} to the homomorphism $G\wr H\to Q\wr H$. By 
hypothesis the FIC$ _{\cal C}$ is true for $Q\wr H$. So let 
$S\in {\cal C}(Q\wr H)$. We have to prove that the 
FIC$ _{\cal C}$ is true for $p^{-1}(S)$. 
Note that $p^{-1}(S)$ contains 
$p^{-1}(S\cap 
Q^H)$ as a normal subgroup of finite index. Therefore using $(2)$ it is enough to 
prove the FICwF$ _{\cal C}$ for $p^{-1}(S\cap
Q^H)$. Next, note that $S\cap
Q^H$ is a subgroup of $\Pi_{h\in H}L_h$, where $L_h$ is the image 
of $S\cap Q^H$ under the projection to the $h$-th coordinate 
of $Q^H$, and since $\cal C$ is closed under taking quotient 
$L_h\in {\cal C}(Q)$. Hence 
$p^{-1}(S\cap Q^H)$ is a subgroup of $\Pi_{h\in H}p^{-1}(L_h)$. 
Since by hypothesis the 
FICwF$ _{\cal C}$ is true for $p^{-1}(L_h)$ 
for $h\in H$, using 
$(1)$, $(2)$ and 
Lemma \ref{finiteindex} we see 
that the FICwF$ _{\cal C}$ is true for $p^{-1}(S\cap
Q^H)$. This completes the 
proof.\end{proof} 

\begin{cor} \label{finabe} ${\cal P} _{\cal 
{VC}}$ implies that the FIC$ _{\cal {VC}}$ is true for 
finitely generated abelian group.\end{cor}

\begin{proof} The proof is immediate from $(1)$ of Proposition 
\ref{product} since the FIC$ _{\cal {VC}}$ is true 
for virtually cyclic groups.\end{proof} 

\begin{rem} {\rm In $(2)$ of Proposition \ref{product} if we assume 
that the FIC$ _{\cal C}(G)$ is satisfied instead of 
the FICwF$_{\cal C}(G)$ then it is not known 
how to deduce the FIC$_{\cal C}(K)$. Even in the case of the 
FIC$_{\cal {VC}}$ and when $G$ is a free group it is open. 
However if 
$G$ is free then the FIC$^{P}(K)$ is satisfied by  
results of Farrell-Jones. See the proof of  
Proposition \ref{graph-finite} for details. Also using a 
recent result of Bartels and L\"{u}ck (\cite{BL1}) it 
can be shown by the same method that the FIC$^L(K)$ is 
satisfied.}\end{rem}

\begin{prop} \label{integer} $_t{\cal T} _{\cal 
C}$ ($_{wt}{\cal T} _{\cal C}$) implies 
that the FIC$ _{\cal C}$ (FICwF$ _{\cal C}$) 
is true for any free group. And 
if $\cal C$ contains the finite groups then 
$_f{\cal T}_{\cal C}$ ($_{wf}{\cal T}_{\cal C}$) 
implies that the FIC$_{\cal C}$ (FICwF$_{\cal C}$) 
is true for a finitely generated group 
which contains a free subgroup of finite index.\end{prop}

\begin{proof} The proof follows from Lemma \ref{finitefree}.\end{proof}

\begin{cor} \label{freeabelian} $_t{\cal T} _{\cal 
C}$ ($_{wt}{\cal T} _{\cal
C}$) and 
${\cal P} _{\cal C}$ imply that the 
FIC$ _{\cal C}$ (FICwF$ _{\cal C}$) is true for 
finitely generated free abelian groups.\end{cor}

\begin{proof} The proof is a combination of Proposition \ref{integer} and 
$(1)$ of Proposition \ref{product}.\end{proof}

Let $F^n$ denote a finitely generated free group of rank $n$.

\begin{lemma} \label{infin} Assume that the FIC$ _{\cal 
C}$ is true for ${\Bbb Z}^{\infty}\rtimes \l t\r$ 
($F^{\infty}\rtimes \l t\r$) for any countable 
infinitely generated abelian group ${\Bbb Z}^{\infty}$, then the 
FIC$ _{\cal C}$ is true for 
${\Bbb Z}^n\rtimes \l t\r$ ($F^n\rtimes \l t\r$) 
for all 
$n\in {\Bbb N}$. Here all actions of $t$ on the corresponding groups 
are arbitrary. 

And the same holds if we replace the FIC$_{\cal C}$ by 
the FICwF$ _{\cal C}$.\end{lemma}

\begin{proof} The proof is an easy consequence of the 
hereditary property and Lemma \ref{finiteindex}. \end{proof}

\begin{lemma} \label{pc} If the FICwF$ _{\cal 
{VC}}$ is true for $F^{\infty}\rtimes \l t\r$ for any 
action of $t$ on $F^{\infty}$, then ${\cal P}_{\cal {VC}}$ is satisfied.\end{lemma}

\begin{proof} Note that ${\Bbb Z}\times {\Bbb Z}$ is a subgroup of
$F^{\infty}\rtimes_t \l
t\r$, where the suffix `t' denotes the trivial action of $\l t\r$ on
$F^{\infty}$. Hence the
FICwF$ _{\cal {VC}}$ is true for 
${\Bbb Z}\times {\Bbb Z}$, that is the FIC$ _{\cal {VC}}$
is true for $({\Bbb Z}\times {\Bbb Z})\wr F$ for any finite group $F$. On
the other hand for two virtually cyclic groups $C_1$ and $C_2$ the
product $C_1\times C_2$ contains a finite 
index free abelian normal subgroup (say $H$) of rank $\leq 2$ (see 
Lemma \ref{last}), and
therefore $C_1\times C_2$ is a subgroup of $H\wr F$ for a finite group
$F$. Using the hereditary property we conclude that 
${\cal P}_{\cal {VC}}$ is satisfied.\end{proof}

\begin{prop} \label{polycyclic} The 
FICwF$^{P}$ is true for any virtually polycyclic group.\end{prop}

\begin{proof} By [\cite{FJ}, Proposition 2.4] the FIC$^{P}$ is true 
for any virtually poly-infinite cyclic group. Also a 
polycyclic group is virtually poly-infinite cyclic. 
Now it is easy to check that the wreath product of a virtually 
polycyclic group with a finite group is virtually poly-infinite cyclic. This 
completes the proof.\end{proof}

Since the product of two virtually cyclic groups is 
virtually polycyclic an immediate corollary is the following. 
But we give a proof of the Corollary independent of 
Proposition \ref{polycyclic}.

\begin{cor}\label{PVC} ${\cal 
P}^{P}$ is satisfied. Also FICwF$^{P}$ is true for 
any virtually cyclic group.\end{cor}

\begin{proof}[Proof of Corollary \ref{PVC}] At first note that 
the FIC$^P$ is true for virtually cyclic groups. Hence for the 
first part we only 
have to prove that the FIC$^P$ is true for $V_1\t V_2$ where 
$V_1$ and $V_2$ are two infinite virtually cyclic groups. 
Note that $V_1\t V_2$ contains a finite index free abelian normal subgroup, 
say $A$, on two generators. Therefore $V_1\t V_2$ embeds in $A\wr H$ for 
some finite group $H$. Since $A$ is isomorphic to the fundamental group 
of a flat $2$-torus, FIC$^P$ is true for $A\wr H$. See Fact 3.1 and Theorem A 
in \cite{FR}. Therefore FIC$^P$ is true for $V_1\t V_2$ by the hereditary 
property. This proves that ${\cal 
P}^{P}$ is satisfied.

The proof of the second part is similar since for any virtually 
cyclic group $V$ and for any finite group $H$, $V\wr H$ 
is either finite or embeds in a group of the type $A\wr H'$ for some finite 
group $H'$ and where $A$ is isomorphic to a free abelian group 
on $|H|$ number of generators and therefore $A$ is isomorphic to the 
fundamental group of a flat $|H|$-torus. Then we can again apply Fact 3.1 
and Theorem A from \cite{FR}.\end{proof}

We will also need the following proposition.

\begin{prop} \label{graph-finite} Let $\cal G$ be a graph of finite
groups. Then $\pi_1({\cal G})$ satisfies the FICwF$^{P}$.\end{prop}

\begin{proof} By Remark \ref{finitegraph} we can assume that the graph is
finite. Lemma \ref{finitefree} implies that we need to show that the
FICwF$^{P}$ is true for finitely generated groups which contains a free
subgroup of finite index. Now it is a formal consequence of results of
Farrell-Jones
that the FICwF$^{P}$ is true for a free
group. For details see [\cite{FL}, Lemma 2.4]. 
Also compare [\cite{FR}, Fact 3.1]. Next since 
${\cal P}^{P}$ is satisfied (Corollary
\ref{PVC}) using $(2)$ of Proposition \ref{product} we complete the
proof. 
\end{proof}

A version of Proposition \ref{graph-finite} can also be proven 
for FICwF$^L$ using the result that the FIC$^L$ is true for a group which 
acts on a finite dimensional $CAT(0)$-space 
from \cite{BL1}. 

\section{Proofs of the Propositions}

Recall that $\cal C$ always denotes a class of groups which is closed 
under isomorphism, taking subgroups and quotients and contains all the
finite groups.

The proofs of the Propositions appear in the following 
sequence: 2.3-2.4-2.1-2.2.

\begin{proof} [Proof of Proposition \ref{maintheorem1}] 
Proof of $(1)$. Since ${\cal H}^?_*$ is continuous by Remark 
\ref{finitegraph} we can assume that the graph $\cal G$ is finite.

$1(a)$. Since the graph $\cal G$ is a 
tree and the vertex groups 
are torsion free by Lemma \ref{stark} the
restriction of the homomorphism 
$f:\pi_1({\cal G})\to H_1(\pi_1({\cal G}), {\Bbb Z})/\{torsion\}=A (say)$ 
to each vertex group is injective. Let $K$ be the kernel 
of $f$. Let $T$ be the tree on which $\pi_1({\cal G})$ acts for the 
graph of group structure $\cal G$. Then $K$ also acts on $T$ 
with vertex stabilizers $K\cap g{\cal G}_v g^{-1}=(1)$ where $g\in 
\pi_1({\cal G})$ and $v\in V_{\cal G}$. Hence by Lemma \ref{finitefree} $K$ is a free 
group (not necessarily finitely generated). Next note that 
$A$ is a finitely generated free abelian group 
and hence the FIC$ _{\cal {VC}}$ is true for 
$A$ by Corollary \ref{finabe} and by hypothesis. Now 
applying Lemma \ref{inverse} to the homomorphism 
$f:\pi_1({\cal G})\to A$ and noting that a torsion 
free virtually cyclic group is either trivial or infinite cyclic (see 
[\cite{FJ95}, Lemma 2.5])  
we complete the proof. We also need to use Lemma \ref{infin} when $K$ 
is finitely generated.

$1(b)$. The proof of this case is almost the same as that of the previous 
case.

Let $A=H_1(\pi_1({\cal G}), {\Bbb Z})$.
Now $A$ is a finitely generated abelian group and hence the 
FIC$ _{\cal {VC}}$ is true for $A$ by 
Corollary \ref{finabe}. Next we apply Lemma \ref{inverse} to 
$p:\pi_1({\cal G})\to A$. Again by Lemma \ref{stark} the kernel $K$ 
of this homomorphism acts on a tree with trivial stabilizers and 
hence $K$ is free. Let $V$ be a virtually cyclic subgroup of 
$A$ with $C<V$ be infinite cyclic subgroup of finite index in 
$V$. Let $C$ be 
generated by $t$. Then the inverse image $p^{-1}(V)$ 
contains $K\rtimes \l t\r$ as a 
normal finite index subgroup. By hypothesis the 
FICwF$ _{\cal {VC}}$ is true for $K\rtimes \l t\r$. 
Now using $(2)$ of Proposition \ref{product} we see that the 
FICwF$ _{\cal {VC}}$ is true for $p^{-1}(V)$ and 
in particular the FIC$ _{\cal {VC}}$ is true for 
$p^{-1}(V)$. This completes the proof of $1(b)$.

$1(c)$. Since the graph $\cal G$ is not a tree it is homotopically 
equivalent to a wedge of $r$ circles for $r\geq 1$. Then there is a surjective 
homomorphism 
$p:\pi_1({\cal G})\to F^r$ where $F^r$ is a free group on $r$ generators. 
And the kernel $K$ of this homomorphism is the fundamental group of the 
universal covering $\widetilde {\cal G}$ of the graph of groups $\cal G$. 
Hence $K$ is the fundamental group of an 
infinite tree of finitely generated abelian groups. Now we would 
like to apply Lemma \ref{inverse} to the homomorphism $p:\pi_1({\cal 
G})\to F^r$. By hypothesis and by Lemma \ref{infin} the FIC$_{\cal {VC}}$ is true 
for any semidirect product $F^n\rtimes \l t\r$, hence the FIC$_{\cal {VC}}$ is 
true for $F^r$ by the hereditary property. Since $F^r$ is torsion free, 
by Lemma \ref{inverse},   
we only have to check that the FIC$ _{\cal {VC}}$  
is true for the semidirect product $K\rtimes \l t\r$ for any action of 
$\l t\r$ on $K$. 

By Lemma \ref{stark} the 
restriction of the 
following map to each vertex group of $\widetilde {\cal G}$ is 
injective. $$\pi_1(\widetilde {\cal G})\to H_1(\pi_1(\widetilde 
{\cal G}), {\Bbb Z}).$$ 

Let us denote the range group and the kernel of the  
homomorphism in the above display by $A$ and $B$ be 
respectively.
Since the commutator subgroup of a group is characteristic we 
have the following exact sequence of groups induced 
by the above homomorphism. $$1\to 
B\to K\rtimes \l t\r \to A\rtimes \l t\r \to 1$$ for any action of $\l 
t\r$ on 
$K$. Recall that $K=\pi_1(\widetilde {\cal G})$. 
Now let $K$ acts on a tree $T$ which induces the tree of 
groups structure $\widetilde {\cal G}$ on $K$. Hence $B$ also acts on $T$ 
with vertex  
stabilizers equal to $B\cap g\widetilde {\cal G}_v g^{-1}=(1)$ 
where $v\in V_{\widetilde {\cal G}}$ and $g\in K$. This follows 
from the fact that the restriction to any vertex group of 
$\widetilde {\cal G}$ of the homomorphism 
$K\to A$ is injective. Thus $B$ acts on 
a tree with trivial stabilizers and hence $B$ is a free group by Lemma 
\ref{finitefree}. 

Next note that the group $A$ is a countable infinitely generated 
abelian group.
Now we can apply Lemma \ref{inverse} to the homomorphism $K\rtimes \l t\r 
\to A\rtimes \l t\r$ (and use Lemma \ref{infin} if $B$ is finitely 
generated) and $(2)$ of Proposition \ref{product} in exactly the same 
way as we did in the proof of $1(b)$. This completes the proof of 
$1(c)$.

{\it Proofs of $2(i)$ and $2(ii)$}. 
Let $e_1,e_2,\ldots , e_k$ be the finite edges of $\cal G$ 
so that each of the connected 
components ${\cal G}_1, {\cal G}_2, \ldots ,{\cal G}_n$ of  
${\cal G}-\{e_1,e_2,\ldots , e_k\}$ is a  
tree of finitely generated abelian groups of the same rank. 
By Lemma \ref{ic} a finite 
tree of finitely generated abelian groups of the same rank  
has the intersection property. Therefore using Lemma \ref{inter} 
we see that such a tree of groups has residually finite fundamental group. 

Now we check that the FICwF$_{\cal {VC}}$ is true for $\pi_1({\cal G}_i)$ 
for $i=1,2,\ldots , n$. 

Assume that ${\cal G}_i$ is a graph of finitely generated abelian 
groups of the same rank (say $r$). Then by Lemma \ref{ic} $\pi_1({\cal G}_i)$ 
contains a finitely generated free abelian normal subgroup $A$ of rank 
$r$ so that the quotient $\pi_1({\cal G}_i)/A$ is isomorphic to the 
fundamental group of a graph of finite-type. Hence  
by the assumption $_{wt}{\cal T} _{\cal {VC}}$ the 
FICwF$ _{\cal {VC}}$ is true for $\pi_1({\cal G}_i)/A$. 
Now we would like to apply $(3)$ of Proposition \ref{product} to the 
following exact sequence. 
$$1\to A\to \pi_1({\cal G}_i)\to \pi_1({\cal G}_i)/A\to 1.$$ 

Note that using $(2)$ of Proposition \ref{product} it is enough 
to prove that the FICwF$ _{\cal {VC}}$ is true for 
$A\rtimes \langle t\rangle$. In case $2(i)$ this follows from 
the hypothesis and Lemma \ref{infin}. And in case $2(ii)$ note that 
$A\rtimes \langle t\rangle$ contains a rank $2$ free abelian 
subgroup of finite index. (See Lemma \ref{last}.) Therefore we can again apply 
$(2)$ of Proposition \ref{product} and Lemma \ref{finabe} 
to see that the FICwF$ _{\cal {VC}}$ is true for
$A\rtimes \langle t\rangle$.

Now observe that there is a finite graph of groups $\cal H$ 
so that the edge groups of $\cal H$ are finite and the 
vertex groups are isomorphic to $\pi_1({\cal G}_i)$ where 
$i$ varies over $1,2,\ldots , n$ and 
$\pi_1({\cal H})\simeq \pi_1({\cal G})$. This follows 
from Lemma \ref{almost} since $\cal G$ is almost a tree of groups.

Next we apply 
$(1)$ of Proposition \ref{residually-prop} 
(since $\pi_1({\cal G}_i)$ is residually finite for $i=1,2,\ldots , n$)  
to $\cal H$ to complete the proof of $(2)$.  

{\it Proof of $(3)$}.
For the proof of $(3)$ of Proposition \ref{maintheorem1} 
recall that in [\cite{S}, Corollary 3.5] it is proved that the 
UNil-groups with respect to free product 
with amalgamation decomposition of $\pi_1({\cal G})$ along 
any edge group of $\cal G$ vanish 
and hence using the 
discussion before Theorem 0.13 in 
\cite{BL} we complete the proof of $(3)$. Since there is 
no UNil-groups we do not need to tensor the assembly map by ${\Bbb 
Z}[\frac{1}{2}]$ in [\cite{BL}, Theorem 0.13].\end{proof}

\begin{proof} [Proof of Proposition \ref{TC}] Lemma \ref{almost} implies that 
there is a graph of groups $\cal H$ with the same fundamental group as 
of $\cal G$ 
so that the edge groups of $\cal H$ are finite and the vertex groups are either 
finite or  
fundamental groups of finite trees of groups with rank $1$ finitely 
generated abelian 
vertex and edge groups. Now note that by Lemma 
\ref{ic} the tree of group (say 
$\cal K$) corresponding to 
a vertex group of $\cal H$ satisfies the hypothesis of Lemma \ref{fi}. 

Therefore $\cal K$ has the intersection property. Also by Lemma 
\ref{rfcyclic} $\pi_1(\cal K)$ is residually finite. Since by 
Corollary \ref{PVC} the FICwF$^{P}$ 
is true for virtually cyclic groups we can
apply $(1)$ and $(3)$ of Proposition \ref{residually-prop} to complete
the proof of the Proposition provided we check that $_{wt}{\cal
T}^{P}$ and ${\cal
P}^{P}$ are satisfied. By Corollary 
\ref{PVC} the last condition is satisfied. 

We now check the first 
condition. So let $\cal G$ be a graph of virtually cyclic 
groups with trivial edge groups. Using Remark \ref{finitegraph} 
we can assume that the graph is finite. Hence by Lemma \ref{almost1} below 
$\pi_1({\cal G})$ is a free product of a finitely generated 
free group and the vertex 
groups of $\cal G$. By Corollary \ref{PVC} 
the FICwF$^{P}$ is true for any virtually 
cyclic group and by Proposition \ref{graph-finite} and 
Lemma \ref{finitefree} it is true for finitely generated free groups. 
Therefore we can apply [\cite{R1}, Theorem 3.1] to conclude that the 
FICwF$^{P}$ is true for $\pi_1({\cal G})$.\end{proof}

Here we remark that when we applied $(3)$ of Proposition \ref{residually-prop} 
in the 
above proof of Proposition \ref{TC} we only needed the fact that 
the FICwF$^P$ be true for the mapping torus of a virtually 
cyclic group (see the proof of $(3)$ of Proposition \ref{residually-prop}). 
We note that for this we do not need Proposition \ref{polycyclic}, instead 
it can easily be deduced from the proof of Corollary \ref{PVC}. We just need 
to mention that the mapping torus of a virtually cyclic group is either 
virtually cyclic or it contains a two generators 
free abelian normal subgroup of finite index (see the following Lemma). 

\begin{lemma} \label{last} Suppose a group $G$ has a virtually cyclic normal subgroup 
$V_1$ with virtually cyclic quotient $V_2$. Then $G$ is virtually cyclic when 
either $V_1$ or $V_2$ is finite and $G$ contains a two-generator free abelian normal 
subgroup of finite index otherwise.\end{lemma}

\begin{proof} We have the following exact sequence. $$1\to V_1\to G\to V_2\to 1.$$ 
For the first assertion there is nothing to prove when $V_2$ is finite.
So let $V_1$ is finite and $V_2$ is infinite and $C$ is an infinite cyclic subgroup of 
$V_2$ of finite index. Let $p$ denote the 
homomorphism $G\to V_2$. Then $p^{-1}(C)$ 
has a finite normal subgroup with quotient $C$. Therefore $p^{-1}(C)$ 
is virtually infinite cyclic. Let $C'$ be an infinite cyclic 
subgroup of $p^{-1}(C)$ 
of finite index. Hence $C'$ is also an infinite cyclic 
subgroup of $G$ of finite 
index. This proves the first assertion. For the second assertion let $V_1$ is also infinite 
and $D$ is an infinite cyclic subgroup of $V_1$ of finite index. Then 
$p^{-1}(C)\simeq D\rtimes C$. Since $C$ and $D$ are both infinite cyclic the 
only possibilities for $p^{-1}(C)$ are that it contains a free abelian subgroup on two 
generators of index either $1$ or $2$. This proves the Lemma.\end{proof}

\begin{proof} [Proof of Proposition \ref{theorem}] Let $T$ be the tree on
which the group $\pi_1({\cal G})$ acts so that the associated graph of
groups is $\cal G$. For the proof, by 
Lemma \ref{inverse} ($(3)$ of Proposition 
\ref{product}) we need to show 
that the FIC$ _{\cal C}$ (FICwF$ _{\cal C}$) 
is true for the homomorphism 
$f$. We prove $(a)$ assuming ${\cal
T} _{\cal C}$ ($_w{\cal
T} _{\cal C}$). 
 
Let $H\in {\cal C}(Q)$. Note that $f^{-1}(H)$ also acts on the tree $T$ with 
stabilizers $f^{-1}(H)\cap\{$stabilizers of the action of $\pi_1({\cal G})\}$. 
Since the restriction of $f$ to the vertex groups  
of ${\cal G}$ have finite kernels we get that 
$f^{-1}(H)\cap\{$stabilizers 
of the action of $\pi_1({\cal G})\}$ is an extension of a finite group 
by a subgroup of $H$. If ${\cal C}={\cal {FIN}}$ then these stabilizers 
also belong to $\cal C$. If ${\cal C}={\cal {VC}}$ then 
using Lemma \ref{last} we see that the stabilizers again belong to 
$\cal C$. 

Therefore the associated graph of groups of the action of $f^{-1}(H)$ 
on the tree $T$ has vertex groups belonging to $\cal C$. 

Now using ${\cal T} _{\cal C}$ ($_w{\cal
T} _{\cal C}$) we 
conclude that 
the FIC$_{\cal C}$ (FICwF$_{\cal C}$) is true for $f^{-1}(H)$. This completes the 
proof of $(a)$.

For the proofs of $(b)$ and $(c)$ just replace  
${\cal T} _{\cal C}$ ($_w{\cal
T} _{\cal C}$) by  
$_f{\cal T} _{\cal C}$ ($_{wf}{\cal
T} _{\cal C}$) and  
$_t{\cal T} _{\cal C}$ ($_{wt}{\cal
T} _{\cal C}$) respectively in the above proof. 
Also use the corresponding assumption on the edge groups of the graph 
of groups $\cal G$.

\end{proof}

\begin{cor} \label{corollary} {\bf (Free products).}
Let $\cal G$ be a finite graph of groups with trivial 
edge groups so that the vertex groups satisfy the FIC$_{\cal C}$ 
(FICwF$_{\cal C}$). If ${\cal P} _{\cal 
C}$ and $_t{\cal T} _{\cal 
C}$ ($_{wt}{\cal T} _{\cal
C}$) are  
satisfied then the FIC$_{\cal C}$ (FICwF$_{\cal C}$) 
is true for $\pi_1({\cal G})$.\end{cor} 

\begin{proof} The proof combines the following two Lemmas.\end{proof}

\begin{lemma} \label{almost1} Let $\cal G$ be a finite graph 
of groups with trivial edge groups. Then there is an
isomorphism
$\pi_1({\cal G})\simeq G_1*\cdots *G_n*F$
where $G_i$'s are vertex groups of ${\cal G}$ 
and $F$ is a free group.\end{lemma}

\begin{proof} The proof is by
induction on the number of edges of the graph. If the graph
has no edge then there is nothing to prove. So assume
${\cal G}$ has
$n$ edges and that the Lemma is true for graphs with $\leq
n-1$ edges.
Let
$e$ be an edge of ${\cal G}$. If ${\cal G}-\{e\}$ is connected
then $\pi_1({\cal G})\simeq \pi_1({\cal
G}_1)*{\Bbb Z}$ where ${\cal G}_1={\cal G}-\{e\}$ is a graph
with $n-1$
edges. On the other hand if ${\cal G}-\{e\}$ has two
components say ${\cal
G}_1$ and ${\cal G}_2$ then $\pi_1({\cal G})\simeq \pi_1({\cal
G}_1)* \pi_1({\cal G}_2)$ where ${\cal G}_1$ and ${\cal G}_2$
has $\leq
n-1$ edges. Therefore by induction we complete the proof of the
Lemma.\end{proof}

\begin{lemma} \label{claim} Assume ${\cal P} _{\cal
C}$ and $_t{\cal T} _{\cal
C}$ ($_{wt}{\cal T} _{\cal
C}$) are
satisfied. If the FIC$ _{\cal C}$ 
(FICwF$ _{\cal C}$) is true 
for  
$G_1$ and $G_2$ then the FIC$ _{\cal C}$ 
(FICwF$ _{\cal C}$) 
is 
true for $G_1 * G_2$.\end{lemma}

\begin{proof} Consider the surjective homomorphism $p:G_1 * G_2\to G_1 
\times G_2$. By $(1)$ of 
Proposition  
\ref{product} the 
FIC$ _{\cal C}$ (FICwF$ _{\cal C}$) 
is true for $G_1\times G_2$. 

Now note that the group $G_1 * G_2$ acts on a tree with 
trivial edge stabilizers and vertex stabilizers conjugates of $G_1$ or 
$G_2$. Therefore the restrictions of $p$ to the stabilizers of this action 
of $G_1*G_2$ on the tree are injective. Hence we are in the situation of 
 $(c)$ of Proposition \ref{theorem},

This completes the proof of the Lemma.\end{proof} 

\begin{rem}{\rm We recall here that in the case of the 
FICwF$^{P}$ the Lemma \ref{claim} coincides with the Reduction Theorem 
(Theorem 3.1) in 
\cite{R1}.}\end{rem}

\begin{proof}[Proof of Proposition \ref{residually-prop}] Since the 
equivariant homology theory is assumed to be continuous when the graph 
is infinite, we can assume
that $\cal G$ is a finite graph of groups.

\noindent
$(1).$ By hypothesis the edge groups of $\cal G$ are finite and the
vertex groups 
are residually finite and satisfy the FICwF$ _{\cal 
{C}}$. By Lemma \ref{graph-res} $\pi_1({\cal G})$ is residually finite. Let 
$F_1,F_2,\ldots ,F_n$ be the edge groups. 
Let $1\neq g\in\cup_{i=1}^nF_i$. Since 
$\pi_1({\cal G})$ is residually finite there is a finite index normal subgroup 
$N_g$ of 
$\pi_1({\cal G})$ so that $g\in \pi_1({\cal G})-N_g$. Let $N=\cap_{g\in 
\cup_{i=1}^nF_i}N_g$. Then $N$ is a finite index normal subgroup of
$\pi_1({\cal G})$ so that $N\cap (\cup_{i=1}^nF_i)=\{1\}$. 

Let $T$ be a tree on which $\pi_1({\cal G})$ acts so that the 
associated graph of group structure on $\pi_1({\cal G})$ is $\cal G$. Hence $N$ 
also acts on $T$. Since $N$ is normal in $\pi_1({\cal G})$ and 
$N\cap (\cup_{i=1}^nF_i)=\{1\}$, the edge stabilizers of 
this action are trivial and the vertex stabilizers are subgroups of
conjugates of 
vertex groups of $\cal G$. Therefore by Corollary \ref{corollary} the FICwF$_{\cal 
{C}}$ is true for $N$. Now $(2)$ of Proposition \ref{product} completes
the proof 
of $(1)$.

\noindent
$(2).$ The proof is by induction on the number of edges. If 
there is no edge then there is one vertex and hence by 
hypothesis the induction starts. Since by hypothesis  
the FICwF$_{\cal {VC}}$ is true for $C\rtimes \l t\r$ for 
$C\in {\cal D}$ it is true for $C$ also 
by Lemma \ref{finiteindex}.
Therefore assume that the result is true for graphs with $\leq n-1$ edges 
which satisfy the hypothesis.
So let $\cal G$ be a finite graph of groups which satisfies the
hypothesis and has $n$ edges. 
Since $\cal G$ has the intersection property there is a normal subgroup
$N$ of $\pi_1({\cal G})$ contained in all the edge groups and of finite
index in 
some edge group, say ${\cal G}_e$ of the edge $e$. 
Let ${\cal G}_1$
be the 
graph of groups with $\cal G$ as the underlying graph and the vertex
and the 
edge groups are ${\cal G}_x/N$ where $x$ is a vertex and an edge  
respectively. Then $\pi_1({\cal G}_1)\simeq \pi_1({\cal G})/N$. 
Let ${\cal G}_2={\cal G}_1-\{e\}$. It is now easy to check that 
the connected components of 
${\cal G}_2$ satisfy the
hypotheses and
also has the intersection property.
Also, since ${\cal G}_2$ has $n-1$ edges, by the induction hypothesis 
$\pi_1({\cal H})$ satisfies the FICwF$_{\cal 
{VC}}$ where $\cal H$ is a connected component of ${\cal G}_2$. 
We now use Lemma \ref{inter} to conclude that $\pi_1({\cal H})$
is also residually finite where $\cal H$ is as above. 
Therefore by $(1)$ $\pi_1({\cal G}_1)$
satisfies the FICwF$ _{\cal 
{VC}}$. Next we apply $(3)$ of Proposition \ref{product} to the homomorphism 
$\pi_1({\cal G})\to \pi_1({\cal G}_1)$.  Note that using $(2)$ of
Proposition \ref{product} it is enough to show the FICwF$_{\cal
{VC}}$ is true for the inverse image of any infinite cyclic 
subgroup of $\pi_1({\cal G}_1)$, but such a group is of the 
form $N\rtimes \l t\r$. Now since $N$ is a subgroup of the edge 
groups of $\cal G$ and $\cal D$ is closed under 
taking subgroups we get $N\in \cal D$, as by hypothesis all the 
vertex groups of $\cal G$ belong to $\cal D$. Again by definition 
of $\cal D$ the FICwF$_{\cal
{VC}}$ is true for $N\rtimes \l t\r$. This 
completes the proof of $(2)$.

\noindent
$(3).$ The proof of $(3)$ follows from $(2)$. For the polycyclic case we
only need to note that virtually polycyclic groups are residually finite
and quotients and subgroups of virtually polycyclic groups are 
virtually polycyclic. Also  
mapping torus of a virtually polycyclic group is again 
virtually polycyclic.

And for the nilpotent case recall from \cite{MR} that the
fundamental group of a finite graph of finitely generated nilpotent
groups is subgroup separable if and only if the graph of groups satisfies
the intersection property. Also note that finitely generated nilpotent
groups are virtually polycyclic.

\noindent
$(4).$ Note that by hypothesis the FICwF$ _{\cal {VC}}$ 
is true for closed surface groups and by 
\cite{GB1} closed surface groups 
are residually finite. Therefore by Lemma \ref{res-surface} 
and by $(1)$ it is enough
to consider a finite graph of groups whose vertex and edge groups are
infinite closed surface groups and the graph of groups satisfies the hypothesis. 
Again by Lemma \ref{res-surface} we have
the exact sequence: $1\to H\to \pi_1({\cal G})\to \pi_1({\cal
G})/H\to 1$. Where $H$ is a closed surface group and $\pi_1({\cal  
G})/H$ is isomorphic to the fundamental group of a finite-type graph of 
groups and hence contains a finitely generated free subgroup of finite
index by Lemma \ref{finitefree}. Hence by $_{wt}{\cal T}_{\cal {VC}}$ and by $(2)$ of Proposition \ref{product}
the FICwF$ _{\cal {VC}}$ is true for $\pi_1({\cal
G})/H$. Now apply $(3)$ of Proposition \ref{product} to the homomorphism
$\pi_1({\cal G})\to \pi_1({\cal  
G})/H$. Using $(2)$ of Proposition \ref{product} it is enough to prove the
FICwF$ _{\cal {VC}}$ for $H\rtimes \l t\r$. Since $H$ is a
closed surface group the action of $\l t\r$ on $H$ is induced by a
diffeomorphism of a surface $S$ so that $\pi_1(S)\simeq H$. Hence
$H\rtimes \l t\r$ is isomorphic to the fundamental group of a closed
$3$-manifold which fibers over the circle. Therefore, using the hypothesis 
we complete the proof.\end{proof}

\section{Proofs of the Theorems}

\begin{proof}[Proof of Theorem \ref{residually}] By Corollary \ref{tt}
 $_{wf}{\cal T}^{P}$ is satisfied and hence $_{wt}{\cal T}^{P}$ 
is also satisfied. Next by Corollary \ref{PVC} ${\cal
P}^{P}$ is satisfied. Therefore the
proof of $(1)$ is completed using $(1)$ of Proposition 
\ref{residually-prop}.

For the proof of $(2)$ apply $(b)$ of Proposition \ref{theorem} 
and Corollary \ref{tt}.
\end{proof}

\begin{proof}[Proof of Theorem \ref{ip}] Corollaries \ref{tt}, \ref{PVC}
and $(2)$ of Proposition \ref{residually-prop} prove the Theorem.
\end{proof}

\begin{proof}[Proof of Theorem \ref{introthm}] Corollaries \ref{tt},
\ref{PVC}, Proposition \ref{polycyclic} and $(3)$ of 
Proposition \ref{residually-prop} prove $(1)$ and
$(2)$ of the Theorem. 

To prove $(3)$ we need only to use $(4)$ of
Proposition \ref{residually-prop}, Corollary \ref{tt} and that the
FICwF$^{P}$ is true for the
fundamental groups of $3$-manifolds fibering over the circle from 
\cite{R1} and \cite{R2}. 

$(4)$ is same as Proposition \ref{TC}.

The proof of $(5)$ follows from $2(i)$ of 
Proposition \ref{maintheorem1}, 
Corollaries \ref{tt} and \ref{PVC}. Since for the 
proof of 
$2(i)$ of 
Proposition \ref{maintheorem1} 
 the weaker hypothesis that `the FICwF$ _{\cal {VC}}$ 
is true for ${\Bbb Z}^n\rtimes \langle t\rangle$ for all $n$' is enough. 
Further 
note that this weaker hypothesis is true for the FICwF$^{P}$ by Proposition \ref{polycyclic}. 
\end{proof}

\section{Examples}

We recall some of the known and unknown results about 
the FIC$^{P}$, FIC$^{L}$, FIC$^{K}$ and the FIC$^{KH}$ 
which are related to this paper. 

\noindent
{\bf Known results:} 
\begin{itemize}
\item $_{wf}{\cal T}^P$ (Corollary \ref{tt}),  
$_{wf}{\cal T}^L$ (\cite{R4}),    
$_{wf}{\cal T}^{K}$ 
([\cite{BL}, Theorem 11.4]) (when the underlying ring 
is regular), ${\cal T}^{KH}$ ([\cite{BL}, Theorem 11.1]) are known.

\item ${\cal P}^{P}$ is satisfied (Corollary \ref{PVC}). It is known that the 
FICwF$^{P}$ is true for any virtually poly-infinite cyclic 
group (Proposition \ref{polycyclic}).  Consequently, the FICwF$^{P}$ is true for 
${\Bbb Z}^n\rtimes \l t\r$ for any $n\in {\Bbb N}$. Also see the footnote in 
Section 1.

\item It is known that if a group $G$ is the limit of a directed 
system of subgroups $\{G_i\}_{i\in I}$ directed under inclusions 
and the FIC$^{X}(G_i)$ is 
satisfied for $i\in I$ then the 
FIC$^{X}(G)$ is also 
satisfied for $X=P,K$ or $L$. See [\cite{FL}, Theorem 7.1]. Compare 
Proposition 
\ref{proposition}. And the same is true for the FIC$^{KH}$. 
See proposition 3.4 and [\cite{BL}, Theorem 11.1].

\item ${\cal {FP}}^X$ is already known 
for $X=P$ (Lemma \ref{claim}) and for $X=K$. ${\cal {FP}}^{KH}$ is also 
known  
(\cite{BL}). And in 
the case of the FIC$^{L}$, it is proved 
recently in \cite{R4}. 
\end{itemize}

\noindent
{\bf Unknown results:}
\begin{itemize}
\item ${\cal T}^{K}$ and ${\cal 
P}^{K}$ are 
not known. See  
[\cite{BL}, Remark 8.5]. Also ${\cal T}^P$ is not known for a general situation. 
Except for some cases which we checked in 
Proposition \ref{TC}.
\item It is not yet proved that the FIC$^{P}$ 
is true for $F^n\rtimes \l t\r$ for $n\in {\Bbb N}\cup \{\infty \}$. 
Though some important special 
cases are known. For example if we assume that the action of $\l t\r$ on 
$F^n$ (here $n$ could be $\infty$) is induced by a diffeomorphism of a
surface then it is proved that 
the FIC$^{P}$
is true for $F^n\rtimes \l t\r$. See \cite{R1} and \cite{R2}. Now for a 
countable infinitely generated torsion free abelian group ${\Bbb 
Z}^{\infty}$ it was shown in [\cite{FL}, Corollary 4.4] that 
the FICwF$^{P}({\Bbb
Z}^{\infty}\rtimes  \l t\r)$ is satisfied if every nearly 
crystallographic group ([\cite{FL}, Definition]) 
satisfies the FIC$^{P}$.\end{itemize}

Below we describe some examples of groups 
for which the results in this article are applicable. Furthermore, 
we show that these groups are new and also they are neither 
$CAT(0)$ nor hyperbolic.

\begin{exm}\label{end}{\rm 
{\bf Graphs of virtually polycyclic groups:} We consider an 
amalgamated free product $H$ of two nontrivial infinite virtually polycyclic groups
over a finite group. Next we recall that in [6] the Fibered 
Isomorphism conjecture in the pseudoisotopy case was
proved for the class of groups
which act cocompactly and properly discontinuously on symmetric 
simply connected nonpositively curved
Riemannian manifolds. In the proof of  
[\cite{BFJP}, Theorem A] it was noted that the condition `symmetric' 
can be replaced by `complete' if we consider torsion free groups. 
See [\cite{FR}, Theorem A] also.
We choose polycyclic groups
so that $H$ does not belong to this class. One example
of such a polycyclic group $S$ is of the type $1\to Z\to S\to Z^2\to 1$ as
described below. Here $Z$ is an infinite cyclic group.

Let G be the Lie group consisting of those
$3 \times 3$ matrices with real number entries whose diagonal entries are all
equal to $1$, entries below the diagonal are all equal to $0$,
and entries above the diagonal are arbitrary. Note that G is
diffeomorphic
to Euclidean 3-space. Let S be the subgroup of G whose entries above the
diagonal are restricted to be integers. Then S is discrete and cocompact.
Let E be the coset space of G by S. It clearly fibers over the 2-torus
with fiber the circle. And the fundamental group of E is S which is
nilpotent
but not abelian. On the other hand in \cite{Y} it was 
shown that the fundamental group of a closed nonpositively curved
manifold which is nilpotent must be abelian. This shows that
$E$ can not support a nonpositively curved Riemannian metric. Now by 
[\cite{KL}, Corollary 2.6] it follows that $S$ can not even  
embed in a  group (called {\it Hadamard groups}) which acts  
discretely and cocompactly on a 
complete simply connected nonpositively curved space (that is a $CAT(0)$-space).
 
Now consider $H_i=S\t F_i$ or $H_i=S\wr F_i$ where $F_i$ 
is a finite group for $i=1,2$.
Next we take amalgamated free product of $H_1$ and $H_2$, $H=H_1*_FH_2$ along some
finite group $F$. Then $H$ does not embed in a Hadamard group as before by 
[\cite{KL}, Corollary 2.6] and $H$ is not virtually 
polycyclic. 
But $H$ satisfies the FICwF$^P$ by $(1)$ of Theorem 1.1 and Proposition 5.4.}\end{exm}

\begin{exm}\label{end0}{\rm {\bf Graphs of residually finite groups with 
finite edge groups:} Let $S$ be the fundamental group of a compact 
Haken $3$-manifold which does not support any nonpositively curved 
Riemannian metric. Such $3$-manifolds can easily be constructed by 
cutting along an incompressible torus in a compact Haken 
$3$-manifold and then gluing differently. See \cite{KL} 
for this kind of construction. Next  
let $H_1$ and $H_2$ be two residually finite 
groups for which the FICwF$^P$ is true and such that $S$ is embedded in 
$H_1$. It is easy to construct such $H_1$, for 
instance take $H_1=S*G*F_1$ or  
$H_1=(S\t G)\wr F_1$ or any such combination where $G$ 
is a finitely generated free group and $F_1$ is a finite group. 
By the same argument as in Example \ref{end} (and 
using \cite{R1} and \cite{R2}) it 
follows that $H=H_1*_FH_2$ (along some finite group $F$) 
satisfies the FICwF$^P$ but is neither virtually polycyclic nor embeds 
in a Hadamard group.}\end{exm}

\begin{exm}\label{end1}{\rm {\bf (Almost) a tree of  
finitely generated
abelian groups where the vertex and the edge groups of any 
component subgraph have the same rank:} Fundamental group of 
such a graph of groups can get very
complicated, for example in the simplest case of amalgamated
free product of two infinite cyclic groups over an infinite
cyclic group, that is $H={\Bbb Z}*_{\Bbb Z} {\Bbb Z}$
produces the $(p,q)$-torus knot group where the two inclusions
${\Bbb Z}\to {\Bbb Z}$ defining the amalgamation are multiplications by $p$ and
$q$, $(p,q)=1$. Though the FICwF$^P$ is
known for knot groups (\cite{R1}), most of the other groups in this class, as 
far as we know, 
are new for which we prove the FICwF$^P$.}\end{exm}

Let us now note that the group $H$ considered in the above examples 
is not hyperbolic as it contains a free abelian subgroup on more 
than one generator.

\begin{rem}{\rm We conclude by remarking  
 that in this paper [\cite{FJ}, Theorem 4.8] (Proposition 
\ref{polycyclic}) is used in the proofs of 
$(1)$, $(2)$ (when the polycyclic or the nilpotent groups are 
not virtually cyclic), $(3)$ and $(5)$ (when the ranks of the abelian groups are $\geq 2$) 
of Theorem \ref{introthm}. See the proof of Corollary \ref{PVC} and the 
discussion after the proof of Proposition \ref{TC}. In this connection we note 
here that using the recent work of Bartels and L\"{u}ck  
in \cite{BL1} all the results in the Introduction of this article can 
be deduced in the $L$-theory case of the Fibered Isomorphism conjecture. 
The same proofs will go through.  
But for this we need to use the $L$-theory version of [\cite{FJ}, Theorem 4.8]  
in the proofs of the particular cases of the items  
of Theorem \ref{introthm} as mentioned above. See \cite{BFL} for 
the proof of [\cite{FJ}, Theorem 4.8] in the $L$-theory case.}\end{rem} 

\newpage
\bibliographystyle{plain}
\ifx\undefined\bysame
\newcommand{\bysame}{\leavevmode\hbox to3em{\hrulefill}\,}
\fi

\end{document}